\documentclass[12pt]{article}
\usepackage{amsmath}
\usepackage{graphicx}
\usepackage{CJK}
\usepackage{multirow}
\usepackage{makecell}

\usepackage{subfigure}
\usepackage[justification=centering]{caption}
\input{amssym.tex}

\def\bc{\begin{center}}       \def\ec{\end{center}}
\def\ba{\begin{array}}        \def\ea{\end{array}}
\def\be{\begin{equation}}     \def\ee{\end{equation}}
\def\bea{\begin{eqnarray}}    \def\eea{\end{eqnarray}}
\def\beaa{\begin{eqnarray*}}  \def\eeaa{\end{eqnarray*}}

\def\mathbb{\Bbb}

\begin{document}
\baselineskip 18pt
\centerline {\bf \large On the number of limit cycles bifurcating from the }
\vskip 0.1 true cm
\centerline {\bf \large  linear center with an algebraic switching curve}

\vskip 0.3 true cm

\centerline{\bf  Jiaxin Wang, Jinping Zhou, Liqin Zhao$^{*}$}
 \centerline{ School of Mathematical Sciences, Beijing Normal University,} \centerline{Laboratory of Mathematics and Complex Systems, Ministry of
Education,} \centerline{Beijing 100875, The People's Republic of China}

\footnotetext[1]{
This work was supported by NSFC(11671040)   \\ * Corresponding author.
E-mail:  zhaoliqin@bnu.edu.cn (L. Zhao).}
\vskip 0.2 true cm

\noindent{\bf Abstract}
This paper studies the family of piecewise linear differential systems in the plane with two pieces separated by a cubic curve. By analyzing the obtained first order Melnikov function, we give an upper bound of the number of limit cycles which bifurcate from the period annulus around the origin under $n$ degree polynomial perturbations. In the case $n=1$ and 2, we obtain that there have exactly 3 and 6 limit cycles bifurcating from the period annulus respectively. The result shows that the switching curves affect the number of the appearing of limit cycles.


\vskip 0.5 true cm
\centerline{\bf{ $\S$1.} Introduction and the main results}
\vskip 0.5 true cm

In the real world, non-smooth phenomena exist in large numbers because of the influence of natural laws and many factors, for more details see, for instance, $\cite{M2008}$ and the references therein.

Piecewise smooth differential system is a kind of important non-smooth system which is based on non-smooth model. Usually, piecewise differential systems have been considered when a straight line separates the plane in two half-planes. In recent years, many authors have studied intensively the number of limit cycles of discontinuous piecewise linear differential systems with two zones separated by a straight line, see for instance [2, 4, 9, 11, 13] and the references quoted in these papers.

We are mainly interested in studying the existence of limit cycles for piecewise linear differential systems with two pieces separated by a nonlinear switching curve. In $\cite{Lb2019}$, the authors considered the family of piecewise linear differential systems in the plane with two pieces separated by a cubic curve. They studied the class of discontinuous piecewise linear differential systems obtained by perturbing up to order 2 in the small parameter $\epsilon$ the linear center $\dot{x}=y$, $\dot{y}=-x$, and they obtained that 7 is a lower bound for the Hilbert number of this family. In $\cite{LZh 2019}$, the authors studied the discontinuous piecewise differential system in the plane with two pieces separated by a curve of $y=x^{n}$.

In the present paper, motivated by the above references, we will study the number of limit cycles for Hamilton system under perturbations of piecewise polynomials of degree $n$ with switching curves $y=x^3$. We consider the following perturbed piecewise smooth differential system
$$
\left(
  \begin{array}{c}
    \dot{x} \\
    \dot{y} \\
  \end{array}
\right)
=\left\{
\begin{aligned}
\left(
  \begin{array}{c}
    y+\epsilon p^+(x,y) \\
    -x+\epsilon q^+(x,y) \\
  \end{array}
\right)
,\ y\geq x^3,\\
\left(
  \begin{array}{c}
    y+\epsilon p^-(x,y) \\
    -x+\epsilon q^-(x,y) \\
  \end{array}
\right)
,\ y<x^3, \\
\end{aligned}
\right.\eqno(1.1)
$$
where $$p^{\pm}(x,y)=\sum^{n}_{i+j=0}a^{\pm}_{i,j}x^{i}y^{j},~~q^{\pm}(x,y)=\sum^{n}_{i+j=0}b^{\pm}_{i,j}x^{i}y^{j}$$ are any polynomials of degree $n$. The first integral of system $(1.1)_{\epsilon=0}$ is
$$
H(x,y)=\frac{1}{2}(x^2+y^2)=\frac{h}{2},~~h\in(0,+\infty).
\eqno(1.2)$$
System $(1.1)_{\epsilon=0}$ has a family of periodic orbits
$$
\begin{aligned}
L_{h}&=\left\{(x,y)|H(x,y)=\frac{h}{2},~h\in(0,+\infty),~y\geq x^3\right\}\\
&\cup\left\{(x,y)|H(x,y)=\frac{h}{2},~h\in(0,+\infty),~y<x^3\right\}\\
:&=L_{h}^{+}\cup L_{h}^{-}.
\end{aligned}
$$

Let $H(n)$ denotes the upper bound of the number of limit cycles bifurcating from the period annulus around the origin for all possible polynomials $p^{\pm}(x,y)$ and $q^{\pm}(x,y)$ up to the first order Melnikov function, taking into account the multiplicity. Our main results are the following theorem.

\vskip 0.2 true cm
\noindent{\bf Theorem 1.1.}
Consider system (1.1), by using the first order Melnikov function in $\epsilon$, we have
$$H(1)=3;~~H(2)=6;~~H(n)\leq 6\left[\frac{n}{2}\right]+6~(n\geq 3).$$

\vskip 0.2 true cm

\noindent{\bf Corollary 1.1.} Let us consider the following perturbed piecewise smooth differential system
$$\left(
  \begin{array}{c}
    \dot{x} \\
    \dot{y} \\
  \end{array}
\right)
=\left\{
\begin{aligned}
\left(
  \begin{array}{c}
    y+\epsilon p^+(x,y) \\
    -x+\epsilon q^+(x,y) \\
  \end{array}
\right)
,\ y\geq x^{1/3},\\
\left(
  \begin{array}{c}
    y+\epsilon p^-(x,y) \\
    -x+\epsilon q^-(x,y) \\
  \end{array}
\right)
,\ y<x^{1/3}, \\
\end{aligned}
\right.$$
where $$p^{\pm}(x,y)=\sum^{n}_{i+j=0}a^{\pm}_{i,j}x^{i}y^{j},~~q^{\pm}(x,y)=\sum^{n}_{i+j=0}b^{\pm}_{i,j}x^{i}y^{j}$$ are any polynomials of degree $n$. Then,
by using the first order Melnikov function in $\epsilon$, we have
$$H(1)=3;~~H(2)=6;~~H(n)\leq 6\left[\frac{n}{2}\right]+6~(n\geq 3).$$

\noindent{\bf Theorem 1.2.}
Let us consider the following perturbed piecewise smooth differential system
$$\left(
  \begin{array}{c}
    \dot{x} \\
    \dot{y} \\
  \end{array}
\right)
=\left\{
\begin{aligned}
\left(
  \begin{array}{c}
    y+\epsilon p^+(x,y) \\
    -x+\epsilon q^+(x,y) \\
  \end{array}
\right)
,\ y\geq x^{m},\\
\left(
  \begin{array}{c}
    y+\epsilon p^-(x,y) \\
    -x+\epsilon q^-(x,y) \\
  \end{array}
\right)
,\ y<x^{m}, \\
\end{aligned}
\right.$$
where $m\in N$, and
$$p^{\pm}(x,y)=\sum^{n}_{i+j=0}a^{\pm}_{i,j}x^{i}y^{j},~~q^{\pm}(x,y)=\sum^{n}_{i+j=0}b^{\pm}_{i,j}x^{i}y^{j}$$ are any polynomials of degree $n$. Then,
by using the first order Melnikov function in $\epsilon$, we prove that there exist $k,l\in N$ such taht $H(n)\le kn+l$.

\noindent{\bf Remark 1.1.}  We conjecture that $H(n)=3n$ in Theorem 1.1.

\vskip 0.2 true cm
This paper is organized as follows. In $\S2$, we will give some preliminaries. In $\S3$, we will prove the Theorem 1.1. First, we obtain the algebraic structure of the first order Melnikov function $M(h)$, which are more complicated than the Melnikov function corresponding to the
discontinuous perturbations with the switching curve is a straight line. Then by direct computation, we get the result for $n\geq3$, and use Chebyshev criterion we get the result for $n=1,~2$.

\vskip 0.5 true cm
\centerline{\bf{ $\S$2}. Preliminaries}
\vskip 0.5 true cm

Consider
$$
\left(
  \begin{array}{c}
    \dot{x} \\
    \dot{y} \\
  \end{array}
\right)
=\left\{
\begin{aligned}
\left(
  \begin{array}{c}
    H_{y}^{+}(x,y)+\epsilon f^+(x,y) \\
    -H_{x}^{+}(x,y)+\epsilon g^+(x,y) \\
  \end{array}
\right)
,\ y\geq \phi(x),\\
\left(
  \begin{array}{c}
    H_{y}^{-}(x,y)+\epsilon f^-(x,y) \\
    -H_{x}^{-}(x,y)+\epsilon g^-(x,y) \\
  \end{array}
\right)
,\ y<\phi(x)\\
\end{aligned}
\right.\eqno(2.1)
$$
where $H^{\pm}$, $f^{\pm}$, $g^{\pm}$ and $\phi(x)$ are all $C^{\infty}$ functions satisfying $\phi(0)=0$, $\epsilon\geq0$ is a small parameter. For system $(2.1)_{\epsilon=0}$ we make the following assumptions:

\noindent{\bf (A1)}: There exists an open interval $J$ such that for each $h\in{J}$, there are two points $A(h)$ and $B(h)$ on the curve $y=\phi(x)$ with $A(h)=(a(h),\phi(a(h)))$, $B(h)=(b(h),\phi(b(h)))$ and satisfying
$$
H^{+}(A(h))=H^{+}(B(h))=h,~~H^{-}(A(h))=H^{-}(B(h)),~~a(h)<0<b(h).
$$

\noindent{\bf (A2)}: There is a family of periodic orbits surrounding the origin with clockwise orientation and denoted by $L_{h}=L_{h}^{+}\cup L_{h}^{-}$, $h\in{J}$ where $L_{h}^{+}$ is defined by $H^{+}(x,y)=h$, $y\geq\phi(x)$ and starting from $A(h)$, ending at $B(h)$, $L_{h}^{-}$ is defined by $H^{-}(x,y)=H^{-}(A(h))$, $y\leq\phi(x)$ and starting from $B(h)$, ending at $A(h)$.

\noindent{\bf (A3)}: Curve $L_{h}^{\pm}$, $h\in{J}$ are not tangent to curve $y=\phi(x)$ at points $A(h)$ and $B(h)$. In other words, for each $h\in{J}$,
$$H_{x}^{\pm}(x,y)+H_{y}^{\pm}(x,y)\phi^{'}(x)\neq0$$
at points $A(h)$ and $B(h)$.

\vskip 0.2 true cm
\noindent{\bf Lemma 2.1.} Under assumptions {\bf (A1)-(A3)}, for the first order Melnikov function of system (2.1), we have
$$
M(h)=\int_{L_{h}^{+}}g^{+}dx-f^{+}dy+\frac{H_{x}^{+}(A)+H_{y}^{+}(A)\phi^{'}(a(h))}{H_{x}^{-}(A)+H_{y}^{-}(A)\phi^{'}(a(h))}\int_{L_{h}^{-}}g^{-}dx-f^{-}dy.
\eqno(2.2)$$

\vskip 0.2 true cm
\noindent{\bf Proof}. Making the change of variables
$$
\left\{
\begin{aligned}
    &x=x,\\
    &z=y-\phi(x),\\
\end{aligned}
\right.
$$
system (2.1) becomes the following $C^{\infty}$ system
$$
\left(
  \begin{array}{c}
    \dot{x} \\
    \dot{z} \\
  \end{array}
\right)
=\left\{
\begin{aligned}
\left(
  \begin{array}{c}
    {\widetilde H}_{z}^{+}(x,z)+\epsilon p^+(x,z) \\
    -{\widetilde H}_{x}^{+}(x,z)+\epsilon q^+(x,z) \\
  \end{array}
\right)
,\ z\geq 0,\\
\left(
  \begin{array}{c}
    {\widetilde H}_{z}^{-}(x,z)+\epsilon p^-(x,z) \\
    -{\widetilde H}_{x}^{-}(x,z)+\epsilon q^-(x,z) \\
  \end{array}
\right)
,\ z<0, \\
\end{aligned}
\right.\eqno(2.3)
$$
where ${\widetilde H}^{\pm}(x,z)=H^{\pm}(x,z+\phi(x))$, ${\widetilde H}^{\pm}_{z}(x,z)=H^{\pm}_{y}(x,z+\phi(x))$ and ${\widetilde H}^{\pm}_{x}(x,z)=H_{x}^{\pm}(x,z+\phi(x))+H_{y}^{\pm}(x,z+\phi(x))\phi^{'}(x)$, and $p^{\pm}(x,z)=f^{\pm}(x,z+\phi(x))$, $q^{\pm}=g^{\pm}(x,z+\phi(x))-\phi^{'}(x)f^{\pm}(x,z+\phi(x))$. It is easy to know system $(2.3)_{\epsilon=0}$ satisfy the following assumptions:

(a) There exist an interval ${\widetilde J}$, and two points ${\widetilde A}(h)=(a(h),0)$ and ${\widetilde B}(h)=(b(h),0)$ such that for $h\in{\widetilde J}$
$$
{\widetilde H}^{+}({\widetilde A}(h))={\widetilde H}^{+}({\widetilde B}(h))={\widetilde h},~~{\widetilde H}^{-}({\widetilde A}(h))={\widetilde H}^{-}({\widetilde B}(h)),~~a(h)<0<b(h).$$

(b) The system $(2.3)_{\epsilon=0}$ has an orbital arc ${\widetilde L}_{h}^{+}$ starting from ${\widetilde A}(h)$ and ending at ${\widetilde B}(h)$ defined by ${\widetilde H}^{+}(x,z)={\widetilde h}$, $z\geq0$, and has an orbital arc ${\widetilde L}_{h}^{-}$ starting from ${\widetilde B}(h)$ and ending at ${\widetilde A}(h)$ defined by ${\widetilde H}^{-}(x,z)=H^{-}({\widetilde A}(h))$, $z<0$.

Therefore, by Theorem 1.1 of \cite{han10}, we can get the first order Melnikov function of system (2.3) is
$$
\begin{aligned}
M(h)&=\int_{{\widetilde L}_{h}^{+}}q^{+}(x,z)dx-p^{+}(x,z)dz+\frac{{\widetilde H}^{+}_{x}({\widetilde A})}{{\widetilde H}^{-}_{x}({\widetilde A})}\int_{{\widetilde L}_{h}^{-}}q^{-}(x,z)dx-p^{-}(x,z)dz\\
&=\int_{L_{h}^{+}}(g^{+}(x,y)-\phi^{'}(x)f^{+}(x,y))dx-f^{+}(x,y)(dy-\phi^{'}(x)dx)\\
&+\frac{H_{x}^{+}(A)+H_{y}^{+}(A)\phi^{'}(a(h))}{H_{x}^{-}(A)+H_{y}^{-}(A)\phi^{'}(a(h))}\int_{L_{h}^{-}}(g^{-}(x,y)-\phi^{'}(x)f^{-}(x,y))dx-f^{-}(x,y)(dy-\phi^{'}(x)dx)\\
&=\int_{L_{h}^{+}}g^{+}dx-f^{+}dy+\frac{H_{x}^{+}(A)+H_{y}^{+}(A)\phi^{'}(a(h))}{H_{x}^{-}(A)+H_{y}^{-}(A)\phi^{'}(a(h))}\int_{L_{h}^{-}}g^{-}dx-f^{-}dy.
\end{aligned}
$$
This ends the proof. $\diamondsuit$

Further, similar to the $\cite{han10}$, if $M(h_{0})=0$ and $M^{'}(h_{0})\neq0$ for some $h_{0}\in(0,+\infty)$, then for $|\epsilon|$ small enough system (2.1) has a unique limit cycle near $L_{h_{0}}$. If $h_{0}$ is a zero of $M(h)$ having an odd multiplicity, then for $|\epsilon|$ small enough system (2.1) has at least one limit cycle near $L_{h_{0}}$. Also, if $M(h)$ has at most $k$ zeros counting multiplicity in $h$ on the interval $(0,+\infty)$, then system (2.1) has at most $k$ limit cycles bifurcating from the annulus $\mathop{\cup}\limits_{h\in(0,+\infty)}L_{h}$.

\vskip 0.2 true cm
\noindent{\bf Definition 2.2.}\cite{Grau11} Let $p_0(x),p_1(x),...,p_{n-1}(x)$ be analytic functions on an open interval $J\subset\mathbb{R}$. The ordered set $(p_0(x),p_1(x),...,p_{n-1}(x))$ is said to be an ECT-system on $J$ if, for all $k=1,2,...,n$, any nontrivial linear combination
$$\alpha_0p_0(x)+\alpha_1p_1(x)+...+\alpha_{k-1}p_{k-1}(x)$$
has at most $k-1$ isolated zeros on $J$ counted with multiplicities.

\vskip 0.2 true cm
\noindent{\bf Lemma 2.3.}\cite{xiong17} The ordered set $(p_0(x),p_1(x),...,p_{n-1}(x))$ is an ECT-system on $J$ if and only if, for each $k=1,2,...,n,$
$$W(p_0,p_1,...,p_{k-1})\neq0,$$ for all $x\in J,$
where $W(p_0,p_1,...,p_{k-1})$ is the Wronskian of functions $p_0(x),p_1(x),...,p_{k-1}(x).$

\vskip 0.2 true cm
\noindent{\bf Lemma 2.4.}\cite{li15} Consider $p+1$ linearly independent analytical functions $f_i:U\rightarrow \mathbb{R}, i=0,1,...,p$, where $U \in \mathbb{R}$ is an interval. Suppose that there exists $j \in {0,,1,...,p}$ such that $f_{j}$ has constant sign. Then there exists $p+1$ constants $C_{i} ,i=0,1,...,p$ such that $f(x)=\sum _{i=0}^{p}C_{i}f_{i}(x)$ has at least $p$ simple zeros in $U$.

\vskip 0.5 true cm
\centerline{\bf{ $\S$3.} Proof of the Theorem 1.1.}
\vskip 0.5 true cm

First, we will obtain the algebraic structure of $M(h)$ of system (1.1). For $h\in(0,+\infty)$ and $i,j\in\mathbb{N}$, we denote
$$
J_{i,j}(h)=\int_{L_{h}^{+}}x^{i}y^{j}dx,~~~I_{i,j}(h)=\int_{L_{h}^{-}}x^{i}y^{j}dx.
$$

\vskip 0.2 true cm
\noindent{\bf Lemma 3.1.} The first order Melnikov function $M(h)$ can be written as
$$
M(h)=\sum^{n}_{i+j=0}\rho_{i,j}^{+}J_{i,j}(h)+\sum^{n}_{i+j=0}\rho_{i,j}^{-}I_{i,j}(h)+\Phi(\sigma(h)),\eqno(3.1)$$
where $\rho_{i,j}^{\pm}$ are arbitrary constants, $\Phi(u)$ is polynomial of $u$ with degree no more than $3n+\frac{5+(-1)^{n}}{2}$.

\vskip 0.2 true cm
\noindent{\bf Proof.} By the Lemma 1.1, we have
$$
\begin{aligned}
M(h)&=\int_{L_{h}^{+}}q^{+}dx-p^{+}dy+\int_{L_{h}^{-}}q^{-}dx-p^{-}dy\\
    &=\sum^{n}_{i+j=0}\left(\int_{L_{h}^{+}}b^{+}_{i,j}x^{i}y^{j}dx-a^{+}_{i,j}x^{i}y^{j}dy+\int_{L_{h}^{-}}b^{-}_{i,j}
    x^{i}y^{j}dx-a^{-}_{i,j}x^{i}y^{j}dy\right).
\end{aligned}
\eqno(3.2)
$$
Notice the curve $y=x^3$ is symmetrical with respect to original point, we suppose that the orbit $L_{h}^{+}(L_{h}^{-})$ intersects the curve $y=x^3$ at points $A(-\sigma(h),-\sigma(h)^{3})$ and $B(\sigma(h),\sigma(h)^{3})$, where
$$\sigma(h)=\frac{\sqrt{6}}{6}\frac{\sqrt{(108h+12\sqrt{81h^2+12})^{1/3}((108h+12\sqrt{81h^2+12})^{2/3}+12)}}{(108h+12\sqrt{81h^2+12})^{1/3}}.$$
Using the Green's Formula, we have
$$
\begin{aligned}
\int_{L_{h}^{+}}x^{i}y^{j}dy&=\int_{L_{h}^{+}\cup\widehat{BOA}}x^{i}y^{j}dy-\int_{\widehat{BOA}}x^{i}y^{j}dy\\
&=-i\iint_{int(L_{h}^{+}\cup\widehat{BOA})}x^{i-1}y^{j}dxdy-\frac{3((-1)^{i+3j+3}-1)}{i+3j+3}\sigma(h)^{i+3j+3},\\
\int_{L_{h}^{+}}x^{i}y^{j}dx&=\int_{L_{h}^{+}\cup\widehat{BOA}}x^{i}y^{j}dx-\int_{\widehat{BOA}}x^{i}y^{j}dx\\
&=j\iint_{int(L_{h}^{+}\cup{\widehat{BOA}})}x^{i}y^{j-1}dxdy-\frac{(-1)^{i+3j+1}-1}{i+3j+1}\sigma(h)^{i+3j+1},
\end{aligned}
$$
which imply that
$$
\int_{L_{h}^{+}}x^{i}y^{j}dy=-\frac{i}{j+1}\int_{L_{h}^{+}}x^{i-1}y^{j+1}dx-\frac{(-1)^{i+3j+3}-1}{j+1}\sigma(h)^{i+3j+3}.\eqno(3.3)
$$
In the similar way, we have
$$
\int_{L_{h}^{-}}x^{i}y^{j}dy=-\frac{i}{j+1}\int_{L_{h}^{-}}x^{i-1}y^{j+1}dx-\frac{1-(-1)^{i+3j+3}}{j+1}\sigma(h)^{i+3j+3}.\eqno(3.4)
$$
From (3.2)--(3.4), we can obtain
$$
\begin{aligned}
M(h)=&\sum^{n}_{i+j=0}\left(\int_{L_{h}^{+}}\left(b^{+}_{i,j}x^{i}y^{j}+\frac{i}{j+1}a^{+}_{i,j}x^{i-1}y^{j+1}\right)dx+\int_{L_{h}^{-}}\left(
b^{-}_{i,j}x^{i}y^{j}+\frac{i}{j+1}a^{-}_{i,j}x^{i-1}y^{j+1}\right)dx\right)\\
&+\sum^{n}_{i+j=0}(a^{+}_{i,j}-a^{-}_{i,j})\frac{(-1)^{i+3j+3}-1}{j+1}\sigma(h)^{i+3j+3}\\
:=&\sum^{n}_{i+j=0}\rho^{+}_{i,j}J_{i,j}(h)+\sum^{n}_{i+j=0}\rho^{-}_{i,j}I_{i,j}(h)+\Phi(\sigma(h)),
\end{aligned}
$$
where $\Phi(\sigma(h))=\sum^{n}_{i+j=0}(a^{+}_{i,j}-a^{-}_{i,j})\frac{(-1)^{i+3j+3}-1}{j+1}\sigma(h)^{i+3j+3}$, $\rho^{\pm}_{i,j}=b^{\pm}_{i,j}+\frac{i+1}{j}a^{\pm}_{i+1,j-1}(j\geq1)$, $\rho^{\pm}_{i,0}=b^{\pm}_{i,0}$. This ends the proof. $\diamondsuit$

\vskip 0.2 true cm
\noindent{\bf Lemma 3.2.}
For $h\in(0,+\infty)$, $i+j=n$ and $l\geq2$, we have

\noindent (i) If $j$ is an even, then $J_{i,j}(h)=-I_{i,j}(h)$; If $j$ is an odd, then $J_{i,j}(h)=I_{i,j}(h)$ for $i$ is an even, $J_{i,j}(h)=-I_{i,j}(h)$ for $i$ is an odd.

\noindent (ii) If $n=2l$, then
$$
J_{i,j}(h)=\alpha_{i,j}h^{[\frac{n}{2}]}J_{0,0}(h)+\beta_{i,j}h^{[\frac{n}{2}]-1}J_{1,1}(h)+\sum^{3[\frac{n}{2}]-3}_{k=0}\varphi_{[\frac{n}{2}]-1-[\frac{k+2}{3}]}(h)\sigma(h)^{7+2k},
\eqno(3.5)$$
\noindent (iii) If $n=2l-1$, then
$$
J_{i,j}(h)=\gamma_{i,j}h^{[\frac{n}{2}]}J_{0,1}(h),
\eqno(3.6)$$
where $\alpha_{i,j}$, $\beta_{i,j}$ and $\gamma_{i,j}$ are arbitrary constants, and $\varphi_{s}(h)$ are polynomials of $h$ with degree no more than $s$.

\vskip 0.2 true cm
\noindent{\bf Proof.}
By direct computation, we can obtain
$$
\begin{aligned}
J_{i,2m}(h)&=\int_{L_{h}^{+}}x^{i}y^{j}dx=\int^{\sigma(h)}_{-\sigma(h)}x^{i}(h-x^{2})^{m}dx=-\int^{-\sigma(h)}_{\sigma(h)}x^{i}(h-x^{2})^{m}dx=-I_{i,2m}(h).
\end{aligned}
$$
It is similar with $J_{2s,2m+1}(h)=I_{2s,2m+1}(h)$ and $J_{2s+1,2m+1}(h)=-I_{2s+1,2m+1}(h)$. Particularly, we have
$$
\begin{aligned}
J_{2s+1,2m}(h)&=\int_{L_{h}^{+}}x^{2s+1}y^{2m}dx=\int_{L_{h}^{+}}x^{2s+1}(h-x^2)^{m}dx=\sum^{m}_{k=0}\frac{m!h^{k}}{k!(m-k)!}\int_{-\sigma(h)}^{\sigma(h)}x^{2(m-k+s) +1}dx\\
&=\sum^{m}_{k=0}\frac{m!h^{k}}{2k!(m-k)!(m-k+s+2)}x^{2(m-k+s+2)}\bigg|_{-\sigma(h)}^{\sigma(h)}=0.
\end{aligned}
$$
It is similar with $I_{2s+1,2m}(h)=0$.
Without loss of generality, we only prove (3.5), and (3.6) can be shown in a similar way. Differentiating $H(x,y)=\frac{h}{2}$ with respect to $x$, we obtain
$$
x+y\frac{\partial y}{\partial x}=0.
\eqno(3.7)$$
Multiplying (1.2) and (3.7) by $x^{i}y^{j}dx$ and $x^{i+1}y^{j}dx$ respectively and integrating over $L_{h}^{+}$, noting (3.3) we have
$$
J_{i+2,j}(h)+J_{i,j+2}(h)=hJ_{i,j}(h)
\eqno(3.8)$$
$$
J_{i+2,j}(h)-\frac{i+1}{j+2}J_{i,j+2}(h)-\frac{(-1)^{i+3j+7}-1}{j+2}\sigma(h)^{i+3j+7}=0.
\eqno(3.9)$$
Elementary manipulations reduce Eps. (3.8) and (3.9) to
$$
J_{i,j}(h)=\frac{i}{i+j+1}\left(hJ_{i,j-2}(h)-\frac{(-1)^{i+3j+1}-1}{j}\sigma(h)^{i+3j+1}\right),\eqno(3.10)$$
$$
J_{i,j}(h)=\frac{j+2}{i+j+1}\left(\frac{i-1}{j+2}hJ_{i-2,j}(h)+\frac{(-1)^{i+3j+5}-1}{j+2}\sigma(h)^{i+3j+5}\right).
\eqno(3.11)$$
We will prove the conclusion by induction on $n$. When $l=2$, (3.10) and (3.11)
give
$$
\begin{aligned}
J_{4,0}(h)&=\frac{1}{5}h^{2}J_{0,0}(h)-\frac{2}{5}h\sigma(h)^{7}-\frac{2}{5}\sigma(h)^{9},\\
J_{3,1}(h)&=\frac{2}{5}hJ_{1,1}(h)-\frac{2}{5}\sigma(h)^{11},\\
J_{2,2}(h)&=\frac{2}{15}h^{2}J_{0,0}(h)-\frac{4}{15}h\sigma(h)^{7}+\frac{2}{5}\sigma(h)^{9},\\
J_{1,3}(h)&=\frac{3}{5}hJ_{1,1}(h)+\frac{2}{5}\sigma(h)^{11},\\
J_{0,4}(h)&=\frac{8}{15}h^{2}J_{0,0}(h)+\frac{8}{15}h\sigma(h)^{7}+\frac{2}{5}\sigma(h)^{13},\\
\end{aligned}
$$
which yield the conclusion for $l=2$. Suppose that the result holds for $l\leq k-1(k\geq4)$. Then for $l=k(k\geq4)$, taking $(i,j)=(0,2k),~(1,2k-1),...,(2k-2,2)$ in (3.10), $(i,j)=(2k-1,1),~(2k,0)$ in (3.11) respectively, we can obtain that
$$
\left(\begin{matrix}
          J_{0,2k}(h)\\
          J_{1,2k-1}(h)\\
          \vdots\\
          J_{2k-2,2}(h)\\
          J_{2k-1,1}(h)\\
          J_{2k,0}(h)
          \end{matrix}\right)\ \
=\left(\begin{matrix}
        \frac{2k}{2k+1}hJ_{0,2k-2}(h)+\frac{2}{2k+1}\sigma(h)^{6k+1}\\
        \frac{2k-1}{2k+1}hJ_{1,2k-3}(h)+\frac{2}{2k+1}\sigma(h)^{6k-1}\\
        \vdots\\
        \frac{2}{2k+1}hJ_{2k-2,0}(h)+\frac{2}{2k+1}\sigma(h)^{2k+5}\\
        \frac{2k-2}{2k+1}hJ_{2k-3,1}(h)-\frac{2}{2k+1}\sigma(h)^{2k+7}\\
        \frac{2k-1}{2k+1}hJ_{2k-2,0}(h)-\frac{2}{2k+1}\sigma(h)^{2k+5}
        \end{matrix}\right).
\eqno(3.12)$$
By inductive hypothesis and (3.12), we have for $i+j=2k$,
$$
J_{i,j}(h)=\alpha_{i,j}h^{k}J_{0,0}(h)+\beta_{i,j}h^{k-1}J_{1,1}(h)+\sum^{3k-3}_{m=0}\varphi_{k-1-[\frac{m+2}{3}]}(h)\sigma(h)^{7+2m},
$$
where $\varphi_{s}(h)$ is a polynomial of $h$ with degree no more than $s$. This ends the proof. $\diamondsuit$

\vskip 0.1 true cm
Substituting (3.5) and (3.6) into (3.1), we get the algebraic structure of the first order Melnikov function $M(h)$.

\vskip 0.2 true cm
\noindent{\bf Lemma 3.3.}
For $h\in(0,+\infty)$, the first order Melnikov function $M(h)$ can be written as
$$
\begin{aligned}
M(h)=&\alpha(h)J_{0,0}(h)+\beta(h)J_{1,1}(h)+\gamma(h)J_{0,1}(h)\\
&+\sum^{n}_{k=2}\sum^{3[\frac{k}{2}]-3}_{m=0}\Psi_{[\frac{k}{2}]-1-[\frac{m+2}{3}]}(h)\sigma(h)^{7+2m}+\Phi(\sigma(h)),
\end{aligned}
\eqno(3.13)$$
where $\alpha(h)$, $\beta(h)$ and $\gamma(h)$ are polynomials of $h$ satisfying ${\rm deg}\alpha(h)$, $\gamma(h)\leq[\frac{n}{2}]$, ${\rm deg}\beta(h)\leq[\frac{n}{2}]-1$, $\Phi(u)$ is a polynomial of $u$ with degree no more than $3n+\frac{5+(-1)^{n}}{2}$, $\Psi_{s}(h)$ is a polynomial of $h$ with degree no more than $s$.

\vskip 0.2 true cm
\noindent{\bf Proof of the Theorem 1.1.}
By direct calculations, we have
$$
\begin{aligned}
J_{0,0}(h)&=2\sigma(h),\\
J_{0,1}(h)&=2\int^{\sqrt{h}}_{0}\sqrt{h-x^2}dx,\\
J_{1,1}(h)&=-\frac{2}{3}(h-\sigma(h)^{2})^{\frac{3}{2}}.
\end{aligned}
\eqno(3.14)$$
Substituting (3.14) into (3.13), and let $x=\sqrt{h}t$, we have
$$
\begin{aligned}
M(h)=&2\alpha(h)\sigma(h)+2h\gamma(h)\int^{1}_{0}\sqrt{1-t^2}dt-\frac{2}{3}\beta(h)(h-\sigma(h)^{2})^{\frac{3}{2}}\\
&+\sum^{n}_{k=2}\sum^{3[\frac{k}{2}]-3}_{m=0}\Psi_{[\frac{k}{2}]-1-[\frac{m+1}{3}]}(h)\sigma(h)^{7+2m}+\Phi(\sigma(h)),
\end{aligned}
\eqno(3.15)$$
where $\int^{1}_{0}\sqrt{1-t^2}dt=\frac{\pi}{4}$.

Let $\sigma(h)=u$ ($i.e.~h=u^2+u^6$), then $M(h)$ in (3.15) can be written as
$$
\begin{aligned}
M(u)=&2u\alpha(u^2+u^6)+\frac{\pi}{2}(u^2+u^6)\gamma(u^2+u^6)-\frac{2}{3}u^9\beta(u^2+u^6)+\Psi_{6[\frac{n}{2}]+1}(u)+\Phi(u),
\end{aligned}
$$
it is easy to check that $M(h)$ and $M(u)$ have the same number of zeros in $(0,+\infty)$. Thus by Lemma 3.3, we can get
$$
H(n)\leq 6\left[\frac{n}{2}\right]+6,~n\geq3.
$$

\vskip 0.2 true cm
If $n=1$, we have
$$M(h)=2(b_{0,0}^{+}-b_{0,0}^-)\sigma(h)+2(b_{0,1}^{+}+a_{1,0}^{+}+b_{0,1}^{-}+a_{1,0}^-)\int _{0}^{\sqrt h} \sqrt{h-x^2}dx+2(a_{0,0}^{-}-a_{0,0}^+)\sigma(h)^3,$$
let $x=\sqrt{h}t$ and $\sigma(h)=u~(i.e. ~h=u^2+u^6)$, then $M(h)$ becomes
$$M(u)=2(b_{0,0}^{+}-b_{0,0}^-)u+\frac \pi 2 (b_{0,1}^{+}+a_{1,0}^{+}+b_{0,1}^{-}+a_{1,0}^-)u^2+2(a_{0,0}^{-}-a_{0,0}^+)u^3+\frac \pi 2 (b_{0,1}^{+}+a_{1,0}^{+}+b_{0,1}^{-}+a_{1,0}^-)u^6.$$

By Lemma 2.3, we compute that for $u>0$,
$$W(u)=u\neq0,$$
$$W(u,u^2)=u^2\neq0,$$
$$W(u,u^2,u^3)=2u^3\neq0,$$
$$W(u,u^2,u^3,u^6)=120u^6\neq0$$
This means that the ordered set $(u,u^2,u^3,u^6)$ is an ECT-system in $u\in (0,+\infty)$. Thus, $M(u)$ has at most three isolated zeros in $(0,+\infty)$.
By Lemma 2.4, we can get that $M(u)$ has at least three isolated zeros in $(0,+\infty)$.

\vskip 0.2 true cm
If $n=2$, we have
$$
\begin{aligned}
M(h)=&2(b_{0,0}^{+}-b_{0,0}^{-})\sigma(h)+2(b_{0,1}^{+}+a_{1,0}^{+}+b_{0,1}^{-}
+a_{1,0}^{-})\int _{0}^{\sqrt h}\sqrt {h-x^2}dx\\
&+(2b_{0,2}^{+}+a_{1,1}^{+}-2b_{0,2}^{-}-a_{1,1}^{-})h\sigma(h)\\
&+(-\frac 23 b_{0,2}^{+}+\frac 23 b_{2,0}^{+}-2a_{0,0}^{+}-\frac 13 a_{1,1}^{+}+\frac 23 b_{0,2}^{-}-\frac 23 b_{2,0}^{-}+2a_{0,0}^{-}+\frac 13 a_{1,1}^{-})\sigma(h)^3\\
&+\frac 23 (-b_{1,1}^{+}-2a_{2,0}^{+}+b_{1,1}^{-}+2a_{2,0}^{-})(h-\sigma(h)^2)^{\frac 32}+2(a_{2,0}^{-}-a_{2,0}^{+})\sigma(h)^5+(a_{1,1}^{-}-a_{1,1}^{+})\sigma(h)^7\\
&+\frac 23 (a_{0,2}^{-}-a_{0,2}^{+})\sigma(h)^9,
\end{aligned}
$$
using the same method above, $M(h)$ can be written as
$$
\begin{aligned}
M(u)=&2(b_{0,0}^{+}-b_{0,0}^{-})u+\frac \pi 2 (b_{0,1}^{+}+a_{1,0}^{+}+b_{0,1}^{-}+a_{1,0}^{-})u^2\\
&+(\frac 43 b_{0,2}^{+}+\frac 23 b_{2,0}^{+}+\frac 23 a_{1,1}^{+}-2a_{0,0}^{+}
-\frac 43 b_{0,2}^{-}-\frac 23 b_{2,0}^{-}-\frac 23 a_{1,1}^{-}+2a_{0,0}^{-})u^3\\
&+2(a_{2,0}^{-}-a_{2,0}^{+})u^5+\frac \pi 2 (b_{0,1}^{+}+a_{1,0}^{+}+b_{0,1}^{-}+a_{1,0}^{-})u^6\\
&+2(b_{0,2}^{+}-b_{0,2}^{-})u^7+\frac 23
(-b_{1,1}^{+}-2a_{2,0}^{+}-a_{0,2}^{+}+b_{1,1}^{-}+2a_{2,0}^{-}+a_{0,2}^{-})u^9.
\end{aligned}
$$
By Lemma 2.3, we compute that for $u>0$,
$$W(u)=u\neq0,$$
$$W(u,u^2)=u^2\neq0,$$
$$W(u,u^2,u^3)=2u^3\neq0,$$
$$W(u,u^2,u^3,u^5)=48u^5\neq0,$$
$$W(u,u^2,u^3,u^5,u^6)=2880u^7\neq0,$$
$$W(u,u^2,u^3,u^5,u^6,u^7)=691200u^9\neq0,$$
$$W(u,u^2,u^3,u^5,u^6,u^7,u^9)=5573836800u^{12}\neq0.$$
This means that the ordered set $(u,u^2,u^3,u^5,u^6,u^7,u^9)$ is an ECT-system in $u\in (0,+\infty)$. Thus, $M(u)$ has at most six isolated zeros in $(0,+\infty)$.
By Lemma 2.4, we can get that $M(u)$ has at least six isolated zeros in $(0,+\infty)$.
That is, $H(2)=6$.

\end{document}